\documentclass{conm-p-l}

\copyrightinfo{2006}{International Press}

\usepackage{graphics}

\newtheorem{theorem}{Theorem}[section]
\newtheorem{lemma}[theorem]{Lemma}

\theoremstyle{definition}
\newtheorem{definition}[theorem]{Definition}

\theoremstyle{remark}
\newtheorem{remark}[theorem]{Remark}

\numberwithin{equation}{section}

\newcommand{\cS}{{\mathcal S}}
\newcommand{\q}{\vec{q}}
\newcommand{\vq}{\vec{q}}

\newcommand{\F}{{\mathcal F}}
\newcommand{\e}{\epsilon}
\newcommand{\vth}{\vartheta}

\newcommand{\A}{{\mathcal A}}

\renewcommand{\k}{\kappa}
\newcommand{\ga}{\gamma}
\newcommand{\hga}{\hat{\gamma}}
\newcommand{\Ga}{\Gamma}

\newcommand{\hH}{\hat{H}}

\newcommand{\dl}{\delta}
\newcommand{\Dl}{\Delta}
\renewcommand{\th}{\theta}

\newcommand{\ra}{\rightarrow}
\newcommand{\al}{\alpha}
\newcommand{\be}{\beta}
\newcommand{\sg}{\sigma}

\newcommand{\pa}{\partial}
\newcommand{\z}{\zeta}

\newcommand{\hq}{\hat{q}}

\newcommand{\La}{\Lambda}

\newcommand{\la}{\lambda}
\newcommand{\bq}{\bar{q}}

\newcommand{\nid}{\noindent}

\newcommand{\bB}{\bar{B}}

\newcommand{\om}{\omega}
\newcommand{\Om}{\Omega}

\newcommand{\lag}{\langle}
\newcommand{\rag}{\rangle}

\renewcommand{\O}{{\mathcal O}}
\newcommand{\hz}{\hat{z}}

\begin{document}

\title[Arnold Diffusion]{Arnold Diffusion of the Discrete Nonlinear Schr\"odinger Equation}

\author{Y. Charles Li}
\address{Department of Mathematics, University of Missouri, 
Columbia, MO 65211}
\curraddr{}
\email{cli@math.missouri.edu}
\thanks{}


\subjclass{37, 34, 35, 78, 76}
\date{}


\dedicatory{Communicated by Y. Charles Li, received April 3, 2006 \\
and, in revised form, June 25, 2006.}

\begin{abstract}
In this article, we prove the existence of Arnold diffusion for an interesting 
specific system -- discrete nonlinear Schr\"odinger equation. The proof is for 
the 5-dimensional case with or without resonance. In higher dimensions, the 
problem is open. Progresses are made by establishing a complete set of 
Melnikov-Arnold integrals in higher and infinite dimensions. The openness lies 
at the concrete computation of these Melnikov-Arnold integrals. New machineries 
introduced here into the topic of Arnold diffusion are the Darboux transformation 
and isospectral theory of integrable systems. 
\end{abstract}

\maketitle








\tableofcontents

\section{Introduction}

For a simple example posed by V. I. Arnold \cite{Arn64}, the existence of 
the so-called Arnold diffusion has been proved (see e.g. \cite{Arn64} 
\cite{AA68} \cite{FM00}). The argument involves two parts: A calculation 
of the Melnikov-Arnold integrals \cite{Arn64} and a transversal intersection 
argument \cite{AA68} supported by a $\la$-lemma \cite{FM00}. Other arguments 
of variational type were also developed \cite{Mat04} \cite{Xia98} \cite{BBB03}
\cite{CY04}. 

Nevertheless, the theory of Arnold diffusion is far from complete \cite{Dou88}
\cite{CG94} \cite{Mar96} \cite{Cre97} \cite{Per98} \cite{BT99} \cite{CV00} 
\cite{DLS00} \cite{FM03} \cite{LW04} \cite{DLS06}. The main challenge is dealing with high 
dimensional specific systems of interest in applications. When the dimensions 
of the KAM tori are large, more Melnikov-Arnold integrals are needed to establish 
Arnold diffusion. Calculating these integrals is a dauting task if not impossible, 
even with computers. In an infinite dimensional phase space, even the dimensions of 
the KAM tori become a challenging issue. Are there infinite dimensional KAM tori ?
in what form ? In the classical setting of a Banach space with angles-momenta 
coordinates, the challenge is how to deal with the perturbation v.s. the decay of 
the sequence of momenta. It is a very interesting problem. 

The aim of the current article is to draw attention to two canonical systems of 
mathematical physics: The discrete nonlinear Schr\"odinger equation (DNLS) and 
its continuous version -- the nonlinear Schr\"odinger equation (NLS). DNLS and 
NLS are integrable systems that describe many different phenomena in physics 
\cite{APT04}. DNLS is an integrable finite difference discretization of NLS. An 
interesting fact about DNLS is that one can choose and change the dimensions of 
the phase space by selecting the number of particles in the discretization.
For a two particle case, under periodic Hamiltonian perturbations, the resulting 
system is 5-dimensional, for which the existence of Arnold diffusion will be 
proved here. For a three particle case, the system is 7-dimensional, and will 
be a good testing ground for a higher dimensional theory.

The integrable theory offers the missing link from low dimension to high dimension
via two powerful and beautiful machineries: Darboux transformation and isospectral 
theory. Darboux transformation generates explicit expressions of separatrices, 
while isospectral theory produces all the Melnikov vectors. Together they provide a 
complete set of Melnikov-Arnold integrals with elegant universal formulae. This is 
the case for both DNLS and NLS \cite{Li92} \cite{Li03} \cite{Li04} \cite{Li99} 
\cite{Li04a}. On the other hand, specific calculation of these integrals is the 
challenge. For the purpose of proving the existence of chaos, often one Melnikov 
integral is enough and easily computable \cite{Li04} \cite{Li04a} \cite{Li04b}
\cite{Li04c} \cite{Li06c}. The reason is that one can utilize locally invariant center 
manifolds instead of KAM tori. For the purpose of Arnold diffusion, local invariance
(i.e. orbits can only enter or leave the submanifold through its boundaries) is 
not enough. This is due to the second part of the argument for Arnold diffusion 
mentioned above. To establish a $\la$-lemma, one needs the tori to be invariant 
(not locally invariant). For a different application of $\la$-lemma in (1). 
establishing shadowing lemma in infinite dimensional autonomous systems, (2).
proving the existence of homoclinic tubes and heteroclinically tubular cycles, 
(3). proving the existence of tubular chaos, (4). proving the existence of chaos 
cascade; we refer the readers to \cite{Li03a} \cite{Li03b} \cite{Li04c} \cite{Li06c}.

The article is organized as follows: Section 2 deals with isospectral theory and 
Darboux transformation for both DNLS and NLS. Section 3 deals with Arnold diffusion 
for a 5-dimensional perturbed DNLS.

\section{Isospectral Theory and Darboux Transformation}

In this section, we are going to present the isospectral theory and Darboux 
transformation for both the discrete nonlinear Schr\"odinger equation (DNLS) and 
the nonlinear Schr\"odinger equation (NLS).

\subsection{Discrete Nonlinear Schr\"odinger Equation}

Consider the discrete nonlinear Schr{\"{o}}dinger equation (DNLS),
\begin{equation}
i \dot{q}_n = {1 \over h^2}[q_{n+1}-2 q_n +q_{n-1}] + |q_n|^2(q_{n+1}+q_{n-1})
-2 \om^2 q_n \ ,
\label{DNLS}
\end{equation}
where $q_n$'s are complex-valued, $i = \sqrt{-1}$ is the imaginary unit, $\om$ is 
a positive parameter, and $q_n$ satisfies the periodic boundary condition and even 
constraint,
\begin{equation}
q_{n+N}=q_n,\ \ \ \ q_{-n}=q_n \ ,  \label{bc}
\end{equation}
where $N$ is a positive integer $N \geq 3$ and $h=1/N$. The DNLS (\ref{DNLS}) is a 
$2(M+1)$-dimensional system, where $M=N/2$ ($N$ even) and $M=(N-1)/2$ ($N$ odd).
The DNLS can be rewritten in the Hamiltonian form
\begin{equation}
i \dot{q}_n = \rho_n \pa H_0/ \pa \overline{q_n} \ , \label{HDNLS}
\end{equation}
where $\rho_n = 1 +h^2|q_n|^2$ and 
\[
H_0={1 \over h^2}\sum^{N-1}_{n=0}
\bigg \{ \overline{q_n}(q_{n+1}+q_{n-1})-{2 \over h^2}(1+\om^2 h^2)\ln \rho_n
\bigg \} \ .
\]
The phase space is defined as 
\[
\cS = \bigg\{ \vq = ({\bf q} , \bar{{\bf q}}) \bigg | \ {\bf q}=(q_0,q_1,...,q_{N-1}), 
\ q_{N-n}=q_n  \ (1 \leq n \leq N-1)  \bigg\} \ . 
\]
In $\cS$ (viewed as a vector space over the real numbers), we 
define the inner product, for any two points $\vq^{\ +}$ and $\vq^{\ -}$,
as follows:
\[
\lag \vq^{\ +}, \vq^{\ -} \rag = \sum_{n=0}^{N-1} ( q_n^{+} 
\overline{q_n^{-}} + \overline{q_n^{+}} q_n^{-}) \ . 
\]
\nid
And the norm of $\vq$ is defined as $\|\vq \ \|^2 = \lag \vq,\vq \rag$. 
\begin{remark}
In the expression of $H_0$, both $\sum^{N-1}_{n=0}[ \overline{q_n}(q_{n+1}+q_{n-1})]$ 
and $I= \frac{1}{h^2}\sum^{N-1}_{n=0}\ln \rho_n$ are constants of motion too. $I$ will be 
used later to establish Arnold diffusion in the non-resonant case. Also the constant 
of motion $D$ given by $D^2=\prod^{N-1}_{n=0}\rho_n$ will play an important role in 
the isospectral theory. In the continuum limit (i.e. $h\ra 0$), the Hamiltonian $H_0$ 
has a limit in the manner: $h H_0 \ra H_c$, where $H_c$ is the Hamiltonian for 
NLS, $H_c=-\int^1_0 [ |q_x|^2+2\om^2|q|^2-|q|^4 ]dx$. Also as $h\ra 0$, $D\ra 1$. 
\end{remark}

\subsection{Isospectral Theory of DNLS}

For more details on the topic of this subsection, see \cite{Li03}. DNLS has the Lax 
pair \cite{AL76}
\begin{eqnarray}
\varphi_{n+1}&=&L_n\varphi_n, \label{DLax1} \\
\dot{\varphi}_n&=&B_n\varphi_n, \label{DLax2}
\end{eqnarray}
\noindent
where
\[
L_n=\left( \begin{array}{cc} z& ihq_n \cr ih\overline{q_n} & 1/z \cr \end{array} \right),
\]
\[
B_n={i\over h^2}\left( \begin{array}{ll}
1-z^2+2i\la h-h^2q_n\overline{q_{n-1}}+\om^2 h^2 & -izhq_n+(1/z)ihq_{n-1}  \cr
-izh\overline{q_{n-1}}+(1/z)ih\overline{q_n} & \frac{1}{z^2}-1+2i\la h+h^2\overline{q_n}
q_{n-1}-\om^2 h^2 \cr \end{array} \right),
\]
and where $z=\exp(i\la h)$. Compatibility of the over determined system 
(\ref{DLax1},\ref{DLax2}) gives the ``Lax representation'' 
\[
\dot{L}_n=B_{n+1}L_n-L_nB_n
\]
of the DNLS (\ref{DNLS}). Focusing our attention upon the discrete spatial part 
(\ref{DLax1}) of the Lax pair, let $M_n = M_n(z,\vq )$ be the $2\times 2$ fundamental 
matrix solution such that $M_0$ is the $2\times 2$ identity matrix. The {\em Floquet 
discriminant} $\Dl : \mathbb{C} \times \cS \ra \mathbb{C}$ is defined by 
\[
\Dl = \frac{1}{D}\ \text{trace } M_N(z,\vq )\ ,
\]
where $D^2=\prod^{N-1}_{n=0}\rho_n$. The isospectral theory starts from the fact 
that for any $z \in \mathbb{C}$, 
$\Dl (z, \vq )$ is a constant of motion of DNLS. An easy way to understand this is that 
as $q_n$ evolves in time according to DNLS, the parameter $z$ in the Lax pair does not 
change. $\Dl (z, \vq )$ is a meromorphic function in $z$ of degree ($+N,-N$), and 
provides $(M+1)$ functionally independent constants of motion, where 
$M = N/2$ ($N$ even), $M=(N-1)/2$ ($N$ odd). There are many ways to generate 
$(M+1)$ functionally independent constants of motion from $\Dl$. The approach which proved to be 
most convenient is by employing all the critical points $z^c$ of $\Dl (z, \vq )$:
\[
\frac{\pa \Dl}{\pa z} (z^c, \q )=0\ .
\]
\begin{definition}
We define a sequence of ($M+1$) constants of motion $F_j$ as follows
\begin{equation}
F_j(\vq ) = \Dl (z^c_j(\vq ), \vq )\ , \quad (j=1, \cdots , M+1)\ .
\label{deff}
\end{equation}
\end{definition}
There is a good description on the locations of these critical points $z^c_j$ in 
the NLS setting \cite{Li04}. These $F_j$'s can be used to build a complete set of 
Melnikov-Arnold integrals for the Arnold diffusion purpose. The Melnikov vectors 
are given by the gradients of these $F_j$'s.
\begin{theorem}
Let $z^c_j(\vq)$ be a simple critical point, then
\begin{equation}
\left (\begin{array}{c} \frac{\pa F_j}{\pa q_n} \cr \cr
\frac{\pa F_j}{\pa \overline{q_n}} \cr \end{array} \right ) = 
{i h (\z -\z^{-1}) \over 
2 W_{n+1}} \left (\begin{array}{c} \psi^{(+,2)}_{n+1}\psi^{(-,2)}_{n}+ 
\psi^{(+,2)}_{n}\psi^{(-,2)}_{n+1} \cr \cr -\psi^{(+,1)}_{n+1}\psi^{(-,1)}_{n}-
\psi^{(+,1)}_{n}\psi^{(-,1)}_{n+1} \cr \end{array} \right )\ , 
\label{fidr}
\end{equation}
where $\psi_n^\pm = (\psi_n^{(\pm,1)}, \psi_n^{(\pm,2)})^T$, ($T$: transpose), are two 
Bloch solutions of the Lax pair (\ref{DLax1},\ref{DLax2}) at ($z^c_j,\vq$) such that 
\[
\psi^\pm_n = D^{n/N} \z^{\pm n/N} \tilde{\psi}^\pm_n\ ,
\]
where $\tilde{\psi}^\pm_n$ are periodic in $n$ with period $N$, $\z$ is a complex constant, 
and $W_n = \ \mbox{det}\ ( \psi^+_n,\psi^-_n )$ is the Wronskian.
\end{theorem}
For a perturbation $H_0+\e H_1$ of the Hamiltonian $H_0$, a complete set of 
Melnikov-Arnold integrals is then given by 
\[
I_j = \int_{-\infty}^{\infty} \sum_{n=0}^{N-1} 2 \ \text{Im} \left [ \frac{\pa F_j}{\pa q_n} \rho_n 
\frac{\pa H_1}{\pa \overline{q_n}} \right ]dt \ , \quad (j=1, \cdots , M+1)\ ,
\]
which are evaluated on the unstable manifolds of tori that persist into KAM tori, 
where $\frac{\pa F_j}{\pa q_n}$ is given by (\ref{fidr}). Expression (\ref{fidr}) is a 
universal expression for all the Melnikov vectors. The challenge is of course how 
to compute these $I_j$'s. The difficulty lies at how to obtain expressions of the orbit $q_n$
and the corresponding $\psi_n^\pm$. By utilizing Darboux transformations, these 
expressions can be obtained even by hand in some cases.  

\subsection{Darboux Transformation of DNLS}

For more details on the topic of this subsection, see \cite{Li03}.
Expressions of unstable manifolds of tori can be generated via Darboux transformations.
For DNLS, such a Darboux transformation was established in \cite{Li92}. 

Let $q_n$ be any solution of DNLS. Pick $z = z^d$ ($|z^d| \neq 1$) at which  
the linear operator $L_n$ in the Lax pair (\ref{DLax1})-(\ref{DLax2}) has two 
linearly independent periodic solutions (or anti-periodic solutions) $\phi^{\pm}_n$ 
in the sense:
\[
\phi^{\pm}_{n+N} = D \phi^{\pm}_n \ , \quad (\text{or } \phi^{\pm}_{n+N} = -D \phi^{\pm}_n) \ ,
\]
where $D$ is defined above. Let 
\[
\phi_n=(\phi_{n1}, \phi_{n2})^T = c^+ \phi_n^+ + c^- \phi_n^-\ ,
\]
where $c^+$ and $c^-$ are complex parameters. Define the matrix $\Ga_n$ by
\[
\Ga_n=\left(\begin{array}{cc} z+(1/z)a_n & b_n \cr c_n &-1/z+z d_n \cr
            \end{array} \right),
\]
\nid
where
\begin{eqnarray*}
a_n &=& {z^d \over \overline{z^d}^2\Dl_n}\bigg [|\phi_{n2}|^2+|z^d|^2|\phi_{n1}|^2
    \bigg ],\\
d_n &=& -{1 \over z^d\Dl_n}\bigg [|\phi_{n2}|^2+|z^d|^2|\phi_{n1}|^2
    \bigg ],\\
b_n &=& {|z^d|^4-1 \over \overline{z^d}^2\Dl_n}\phi_{n1} \overline{\phi_{n2}}, \\
c_n &=& {|z^d|^4-1 \over |z^d|^2\Dl_n}\overline{\phi_{n1}}\phi_{n2}, \\
\Dl_n &=& -{1 \over \overline{z^d}}\bigg [|\phi_{n1}|^2+|z^d|^2|\phi_{n2}|^2
    \bigg ].
\end{eqnarray*}
\nid
From these formulae, we see that
\[
\overline{a_n}=-d_n,\ \ \overline{b_n}=c_n.
\]
\begin{theorem}[Darboux Transformation \cite{Li92}]
Define $Q_n$ and $\Psi_n$ by
\[
Q_n = {i\over h}b_{n+1}-a_{n+1}q_n\ , \quad \Psi_n = \Ga_n \psi_n \ ,
\]
where $\psi_n$ solves the Lax pair (\ref{DLax1})-(\ref{DLax2}) 
at $(q_n,z)$. Then $Q_n$ is also a solution of DNLS, and $\Psi_n$ 
solves the Lax pair (\ref{DLax1})-(\ref{DLax2}) at $(Q_n,z)$.
\end{theorem}
In principle, by choosing $q_n$ to be orbits on the tori, one can 
generate $Q_n$ to be orbits on the unstable manifolds of the tori.
In special cases, $Q_n$ can be calculated by hands.

Next we present an example. Define the 2-dimensional invariant plane
\begin{equation}
\Pi = \left \{ \vq \in \cS \ | \ q_n =q\ , \quad \forall n \right \}\ .
\label{invp}
\end{equation}
On $\Pi$, the solutions of DNLS are given by the periodic orbits (1-tori)
\begin{equation}
q_n=q_c \ ,\ \forall n; \quad q_c=a\exp \left \{-i[2(a^2-\om^2)t - \ga]
\right \}\ , 
\label{us}
\end{equation}
where $a$ and $\ga$ are real constants. We choose the amplitude $a$ in 
the following range so that the unstable direction of $q_c$ is 1-dimensional
\begin{eqnarray}
& & N\tan{\pi \over N}< a <N\tan{2\pi \over N}\ ,\ \ \ \mbox{when}\ N>3 \ ,
\nonumber \\ \label{constr} \\
& & 3\tan{\pi \over 3}< a < \infty\ ,\ \ \ \mbox{when}\ N=3 \ .\nonumber
\end{eqnarray}
Increasing the unstable dimensions of $q_c$ amounts to iterations of 
the Darboux transformation which are still doable by hands, and does not 
add difficulty substantially in the Arnold diffusion problem. To apply the 
Darboux transformation we choose 
\[
z^d = \sqrt{\rho}\cos \frac{\pi}{N} +\sqrt{\rho \cos^2 \frac{\pi}{N} -1}\ ,
\]
where $\rho = 1 +h^2a^2$. This $z^d$ is also a critical point of $\Dl$. We 
label it by $z^c_1 = z^d$. Direct calculation leads to the following 
formulae
\begin{eqnarray}
Q_n &=& q_c \bigg [ \Ga / \La_n -1 \bigg ]\ , \label{horbit} \\ \nonumber \\
\left (\begin{array}{c} \frac{\pa F_1}{\pa q_n} \cr \cr
\frac{\pa F_1}{\pa \overline{q_n}} \cr \end{array} \right )\bigg |_{Q_n}
    &=& K [K_n]^{-1} \ \mbox{sech}[2\mu t + 2p] \left ( \begin{array}{c} 
X^1_n \cr X^2_n \cr \end{array} \right ) \ , \label{melv}
\end{eqnarray}
where
\begin{eqnarray*}
\Ga &=& 1-\cos 2 \varphi - i \sin 2 \varphi \tanh [ 2 \mu t + 2p]\ , \\
\La_n &=& 1 \pm \cos \varphi [\cos \be ]^{-1}\ \mbox{sech}[2 \mu t + 2p] 
\cos [2n\be]\ , \\
K &=& -2N (1-\hz^4) [8a\rho^{3/2} \hz^2]^{-1} \sqrt{\rho \cos^2\be - 1}\ , \\
K_n &=& \bigg [ \cos \be \pm \cos \varphi \ \mbox{sech}[2\mu t +2p] 
\cos[2(n-1)\be]\bigg ] \times \\
& &\bigg [ \cos \be \pm \cos \varphi \ \mbox{sech}[2\mu t +2p] 
\cos[2(n+1)\be]\bigg ]\ , \\
X^1_n &=& \bigg [ \cos \be \ \mbox{sech}\ [ 2\mu t +2p] 
\pm (\cos \varphi \\
& & -i \sin \varphi \tanh [ 2\mu t +2p]) \cos [2n\be]\bigg ] e^{i2\th(t)}\ , \\
X^2_n &=& \bigg [ \cos \be \ \mbox{sech}\ [ 2\mu t +2p] 
\pm (\cos \varphi \\
& & +i \sin \varphi \tanh [ 2\mu t +2p]) \cos [2n\be]\bigg ] e^{-i2\th(t)}\ , \\
& & \be = \pi / N\ , \ \ \mu = 2h^{-2} \sqrt{\rho}\sin \be \sqrt{\rho \cos^2 \be -1}\ ,\\
& & \rho = 1+h^2 a^2\ , \ \ h=1/N, \ \ \hz =\sqrt{\rho}\cos \be +\sqrt{\rho \cos^2 \be -1}\ , \\
& & \th(t)=(a^2-\om^2)t - \ga/2\ , \ \ h a e^{i \varphi} =  
\sqrt{\rho \cos^2 \be -1} + i \sqrt{\rho} \sin \be \ . 
\end{eqnarray*}
and $p$ is a real parameter. One can easily see that $Q_n$ represents homoclinic 
orbits asymptotic to the periodic orbits $q_c$: As $t \ra \pm \infty$,
\[
Q_n \ra q_c e^{i(\pi \pm 2\varphi )} \ .
\]
The union 
\[
\bigcup_{\ga \in [0, 2\pi ]}Q_n 
\]
represents the 2-dimensional unstable (=stable) manifold of the 1-torus (\ref{us}). 
When $N>3$, the 1-torus also has a center manifold of codimension 2.

\subsection{Nonlinear Schr\"odinger Equation}

Consider the nonlinear Schr{\"{o}}dinger equation (NLS),
\begin{equation}
iq_t = q_{xx} + 2[|q|^2 - \om^2] q \ ,
\label{NLS}
\end{equation}
where $q = q(t,x)$ is a complex-valued function of the two real 
variables $t$ and $x$, $t$ represents time, and $x$ represents
space. $q(t,x)$ is subject to periodic boundary condition of period 
$2 \pi$, and even constraint, i.e. 
\[
q(t,x + 2 \pi) = q(t,x)\ , \ \ q(t,-x) = q(t,x)\ .
\]
$\om$ is a positive parameter. The DNLS (\ref{DNLS}) is an integrable finite 
difference discretization of the NLS. The NLS can be rewritten in the Hamiltonian form
\begin{equation}
i \dot{q} = \pa H_0/ \pa \bq \ , \label{HNLS}
\end{equation}
where 
\[
H_0=-\int^{2\pi}_0 \left [ |q_x|^2+2\om^2|q|^2-|q|^4 \right ] dx \ .
\]
The phase space is defined as 
\[
\cS = \bigg\{ \vq = (q,\bq ) \bigg | \ q \in H^k \ , \quad (k \geq 1) \bigg\} \ ,
\]
where $H^k$ is the Sobolev space of periodic and even functions.

\subsection{Isospectral Theory of NLS}

For more details on the topic of this subsection, see \cite{Li04}. NLS has the Lax 
pair
\begin{eqnarray}
\psi_x &=& U \psi \ , \label{ZS1} \\
\psi_t &=& V \psi \ , \label{ZS2}
\end{eqnarray}
where
\[
U = i \left ( \begin{array}{lr} \la & q \cr \bq & -\la \cr 
\end{array} \right ) \ , 
\]
\[
V = i \left ( \begin{array}{lr} 2\la^2 -|q|^2 +\om^2 & 
2\la q -i q_x \cr 2 \la \bq + i \overline{q_x} & -2 \la^2 +|q|^2 -\om^2 \cr 
\end{array} \right ) \ . 
\]
Focusing our attention on the spatial part (\ref{ZS1}) of the Lax 
pair (\ref{ZS1},\ref{ZS2}), we can define the fundamental matrix 
solution $M(x)$ such that $M(0)$ is the $2\times 2$ identity matrix. Then 
the {\em Floquet discriminant} $\Dl$ is defined as
\[
\Dl = \ \mbox{trace}\ M(2\pi)\ .
\]
The isospectral theory starts from the fact that for any $\la \in \mathbb{C}$, 
$\Dl (\la , \vq )$ is a constant of motion of NLS. $\Dl = \Dl (\la, q)$ is 
an entire function in both $\la$ and $\vq$. $\Dl$ provides enough functionally 
independent constants of motion to make the NLS (\ref{NLS}) integrable in the 
classical Liouville sense. There are many ways to generate these
functionally independent constants of motion from $\Dl$. The approach which proved to be 
most convenient is by employing all the critical points $\la^c$ of $\Dl (\la ,\vq )$:
\[
\frac{\pa \Dl}{\pa \la} (\la^c, \vq )=0\ .
\]
For each $\vq$, there is a sequence of critical points $\{ \la^c_j \}$ whose 
locations can be estimated \cite{Li04}. 
\begin{definition}
We define a sequence of constants of motion $F_j$ as follows
\[
F_j(\vq ) = \Dl (\la^c_j(\vq ), \vq )\ .
\]
\end{definition}
These $F_j$'s can be used to build a complete set of Melnikov-Arnold integrals for 
the Arnold diffusion purpose. The Melnikov vectors are given by the gradients of 
these $F_j$'s.
\begin{theorem}
Let $\la^c_j(\vq)$ be a simple critical point, then
\begin{equation}
\left (\begin{array}{c} \frac{\pa F_j}{\pa q} \cr \cr
\frac{\pa F_j}{\pa \overline{q}} \cr \end{array} \right ) = 
i \frac {\sqrt{\Dl^2 -4}} {W(\psi^+, \psi^-)}
\left ( \begin{array}{c} \psi_2^+ \psi_2^- \cr - \psi_1^+ \psi_1^-
\cr \end{array} \right ) \ ,
\label{nlsv}
\end{equation}
where $\psi^{\pm} = (\psi^{\pm}_1, \psi^{\pm}_2)^T$, ($T$: transpose), are two 
Bloch solutions of the Lax pair (\ref{ZS1})-(\ref{ZS2}) at ($\la^c_j, \vq$) such that 
\[
\psi^{\pm}(x) = e^{\pm \sg x} \tilde{\psi}^{\pm}(x)\ ,
\]
where $\sg$ is a complex constant and $\tilde{\psi}^{\pm}(x)$ are periodic in $x$ 
with period $2\pi$, 
\[
W(\psi^+,\psi^-)=\psi^+_1 \psi^-_2 - \psi^+_2 \psi^-_1
\]
is the Wronskian, and $\Dl$ is evaluated at $\la = \la^c_j$. 
\end{theorem}
For a perturbation $H_0+\e H_1$ of the Hamiltonian $H_0$, a complete set of 
Melnikov-Arnold integrals is then given by 
\[
I_j = \int_{-\infty}^{\infty} \int_0^{2\pi} 2 \ \text{Im} \left [ \frac{\pa F_j}{\pa q}  
\frac{\pa H_1}{\pa \overline{q}} \right ]dxdt \ , 
\]
which are evaluated on the unstable manifolds of tori that persist into KAM tori, 
where $\frac{\pa F_j}{\pa q}$ is given by (\ref{nlsv}). Expression (\ref{nlsv}) is a 
universal expression for all the Melnikov vectors. The challenge is of course how 
to compute these $I_j$'s. The difficulty lies at how to obtain expressions of the orbit $q$
and the corresponding $\psi^\pm$. By utilizing Darboux transformations, these 
expressions can be obtained even by hand in some cases. Also as mentioned in the 
introduction, there is the issue of elusiveness of infinite dimensional KAM tori. 
As $j \ra \infty$, The magnitude of $I_j$ v.s. the size of the perturbation $\e$ 
is another tricky problem.  

\subsection{Darboux Transformation of NLS}

For more details on the topic of this subsection, see \cite{Li04}.
Expressions of unstable manifolds of tori can be generated via Darboux transformations.
For NLS, such a Darboux transformation was also known. 

Let $q$ be any solution of NLS. Pick $\la = \nu$ (complex constant) at which  
(\ref{ZS1}) has two linearly independent periodic solutions (or anti-periodic solutions) 
$\phi^{\pm}$ in the sense:
\[
\phi^{\pm}(x+2\pi ) = \phi^{\pm}(x)\ , \quad (\text{or } \phi^{\pm} (x+2\pi )= - 
\phi^{\pm}(x)) \ .
\]
Let 
\[
\phi =c_+ \phi^+ + c_- \phi^-\ ,
\]
where $c_+$ and $c_-$ are complex parameters. Define the matrix $G$ by
\[
G = \Ga \left ( \begin{array}{lr} \la -\nu & 0 \cr 
0 & \la - \bar{\nu} \cr \end{array} \right )  \Ga^{-1} \ ,
\]
where 
\[
\Ga = \left ( \begin{array}{lr} \phi_1 & - \overline{\phi_2} \cr 
\phi_2 & \overline{\phi_1} \cr  \end{array} \right ) \ .
\]
\begin{theorem}[Darboux Transformation]
Define $Q$ and $\Psi$ by
\[
Q = q + 2 (\nu - \bar{\nu}) \frac {\phi_1 \overline{\phi_2}} {|\phi_1|^2 
+ |\phi_2|^2}\ , \quad \Psi = G \psi \ ,
\]
where $\psi$ solves the Lax pair (\ref{ZS1})-(\ref{ZS2}) 
at $(q,\la )$. Then $Q$ is also a solution of NLS, and $\Psi$ 
solves the Lax pair (\ref{ZS1})-(\ref{ZS2}) at $(Q,\la )$.
\end{theorem}
In principle, by choosing $q$ to be orbits on the tori, one can 
generate $Q$ to be orbits on the unstable manifolds of the tori.
In special cases, $Q$ can be calculated by hands.

Next we present an example. Define the 2-dimensional invariant plane
\[
\Pi = \left \{ \vq \in \cS \ | \ \pa_x q = 0 \right \}\ .
\]
On $\Pi$, the solutions of NLS are given by the same periodic orbits (1-tori)
as in (\ref{us})
\begin{equation}
q_c=a\exp \left \{-i[2(a^2-\om^2)t - \ga]\right \}\ , 
\label{usc}
\end{equation}
where $a$ and $\ga$ are real constants. We choose the amplitude $a$ in 
the range $a \in (1/2, 1)$ so that the unstable direction of $q_c$ is 1-dimensional.
Increasing the unstable dimensions of $q_c$ amounts to iterations of 
the Darboux transformation which are still doable by hands \cite{Li04a}, and does not 
add difficulty substantially in the Arnold diffusion problem. To apply the 
Darboux transformation we choose 
\[
\nu = i \sg \ , \quad \sg = \sqrt{a^2 - 1/4}\ .
\]
This $\nu$ is also a critical point of $\Dl$. We label it by $z^c_1 = \nu$.
Direct calculation leads to
\begin{eqnarray*}
Q &=& q_c \bigg [ 1 \pm \sin \vth_0 \ \mbox{sech} \tau  
\cos x \bigg ]^{-1} \\
& & \cdot \bigg [ \cos 2\vth_0 - i \sin 2\vth_0 \tanh \tau \mp 
\sin \vth_0 \ \mbox{sech} \tau \cos x \bigg ] \ ,  \\
\left (\begin{array}{c} \frac{\pa F_1}{\pa q} \cr \cr
\frac{\pa F_1}{\pa \overline{q}} \cr \end{array} \right )\bigg |_{Q} 
&=& \frac{1}{4} a^{-2} i (\nu - \bar{\nu})
\sqrt{\Dl(\nu)\Dl''(\nu)} \frac {1}{(|u_1|^2 + |u_2|^2)^2}
\left ( \begin{array}{c} \overline{q_c}\ \overline{u_1}^{\ 2} \cr 
- q_c\ \overline{u_2}^{\ 2} \cr \end{array} \right ) \ ,
\end{eqnarray*}
where 
\begin{eqnarray*}
u_1 &=&  \cosh \frac{\tau}{2} \cos z -i \sinh \frac{\tau}{2} \sin z \ , \\
u_2 &=&  -\sinh \frac{\tau}{2} \cos (z - \vth_0) + i \cosh  \frac{\tau}{2}
\sin (z - \vth_0) \ , \\
& & a e^{i\vth_0} =\frac{1}{2} + \nu \ , \quad \tau = 2 \sg t - \rho \ , 
\quad z = x/2 + \vth_0/2 \mp \pi /4 \ ,
\end{eqnarray*}
and $\rho$ is a real parameter. As $t \ra \pm \infty$,
\[
Q \ra q_c e^{\mp i2\vth_0} \ .
\]
The union 
\[
\bigcup_{\ga \in [0, 2\pi ]}Q 
\]
represents the 2-dimensional unstable (=stable) manifold of the 1-torus (\ref{usc}). 
The 1-torus also has a center manifold of codimension 2. 

\section{Arnold Diffusion}

To establish the existence of Arnold diffusion, one needs three ingredients:
(1). Melnikov-Arnold integrals, (2). A $\la$-lemma, (3). A transversal intersection 
argument. Melnikov-Arnold integrals have been studied above. Next we discuss the 
other two ingredients.
\begin{lemma}[The $\la$-lemma of Fontich-Martin \cite{FM00}]
Let $F$ be a $C^2$ diffeomorphism in $\mathbb{R}^n$, $\mathbb{T}$ be a $C^2$ 
invariant torus in $\mathbb{R}^n$. The dynamics on $\mathbb{T}$ is quasi-periodic.
$\mathbb{T}$ has $C^2$ unstable, stable and center manifolds ($W^u,W^s,W^c$). Let 
$\Ga$ be a $C^1$ manifold intersecting transversally $W^s$ at a point. Then
\[
W^u \subset \overline{\bigcup_{m \geq 0}F^m(\Ga )}\ .
\]
\label{lal}
\end{lemma}
\begin{remark}
Like every other $\la$-lemma, the claim is very intuitive, but the proof is always 
delicate. The proof in \cite{FM00} takes about eight pages. For a quick glance of the 
basic idea, see e.g. \cite{Li03a}. It turns out to be crucial to use Fenichel's 
fiber coordinates. With respect to base points, Fenichel fibers drop one degree of 
smoothness. That is why $C^2$ smoothness is required in the lemma to obtain $C^1$ 
families of $C^2$ unstable and stable Fenichel fibers. Using the Fenichel's 
fiber coordinates, $W^u$ and $W^s$ are rectified, i.e. they coincide with their 
tangent bundles. This makes the estimate a lot easier. Fenichel fibers also make it 
easier to track orbits inside $W^u$ and $W^s$ via the fiber base points in $\mathbb{T}$.
The main argument is to track the tangent space of a submanifold $S$ of $\Ga$ 
starting from the intersection point of $\Ga$ and $W^s$. After enough iterations of $F$, 
one can obtain some estimate of closeness to $W^u$. The novelty of \cite{FM00} is that 
they also track the tangent space at every point in $S$ and off $W^s$. This is 
necessary in order to obtain the claim of the lemma. The claim is proved by showing 
that any neighborhood of any point on $W^u$ has a nonempty intersection with $F^m(S)$ 
for some $m$. The claim of the lemma should also be true in proper infinite 
dimensional settings. 
\end{remark}
\begin{definition}[Transition Chain]
A finite or infinite sequence of tori $\{ \mathbb{T}_j \}$ forms a transition chain if 
$W^s(\mathbb{T}_j)$ intersects transversally $W^u(\mathbb{T}_{j+1})$ at some point, for 
all $j$, and dynamics on $\mathbb{T}_j$ is quasi-periodic.
\end{definition}
\begin{lemma}[Arnold \cite{AA68}]
Let $\{ \mathbb{T}_j \}$ ($1 \leq j \leq N$) be a finite transition chain. Then an 
arbitrary neighborhood of an arbitrary point in $W^u(\mathbb{T}_1)$ is connected to an 
arbitrary neighborhood of an arbitrary point in $W^s(\mathbb{T}_N)$ by an orbit.
\label{al}
\end{lemma}
\begin{proof}
Let $\Om_1$ be an arbitrary neighborhood of an arbitrary point $u_1 \in W^u(\mathbb{T}_1)$.
Let $\bB_{u_1}(r_1) \subset \Om_1$ be a closed ball of radius $r_1>0$ centered at $u_1$. 
Then using the $\la$-lemma \ref{lal}, one can find 
\begin{equation}
\bB_{u_N}(r_N) \subset \bB_{u_{N-1}}(r_{N-1}) \cdots \subset \bB_{u_2}(r_2)
\subset \bB_{u_1}(r_1)
\label{bnt}
\end{equation}
such that
\[
u_j \in W^u(\mathbb{T}_j)\ , \quad r_j >0\ , \quad (1 \leq j \leq N)\ .
\]
Indeed, by the $\la$-lemma \ref{lal},
\[
B_{u_1}(r_1) \cap W^u(\mathbb{T}_2) \neq \emptyset
\]
where $B_{u_1}(r_1)$ is the open ball. Then one can find 
\[
u_2 \in B_{u_1}(r_1) \cap W^u(\mathbb{T}_2) \ , \text{ and }
\bB_{u_2}(r_2) \subset B_{u_1}(r_1)\ .
\]
Again by the $\la$-lemma \ref{lal},
\[
B_{u_2}(r_2) \cap W^u(\mathbb{T}_3) \neq \emptyset \ ,
\]
and one can find 
\[
u_3 \in B_{u_2}(r_2) \cap W^u(\mathbb{T}_3) \ , \text{ and }
\bB_{u_3}(r_3) \subset B_{u_2}(r_2)\ ,
\]
and so on. Let $\Om_N$ be an arbitrary neighborhood of an arbitrary point in 
$W^s(\mathbb{T}_N)$. Inside $\Om_N$, one can find a $C^1$ submanifold intersecting 
transversally $W^s(\mathbb{T}_N)$ at the point. Again by the $\la$-lemma \ref{lal},
\[
B_{u_N}(r_N) \cap F^m(\Om_N) \neq \emptyset \text{  for some } m \ .
\]
Thus $\Om_1$ and $\Om_N$ are connected by an orbit.
\end{proof}
\begin{remark}
It is easy to see that around the connecting orbit, there is in fact a connecting flow 
tube \cite{Li06} \cite{Li03b} \cite{Li04c} \cite{Li06c}
\[
\bigcup_{0\leq m \leq M}F^m(D) \quad \text{for some } M
\]
where $D$ is a neighborhood. When $N=\infty$, relation (\ref{bnt}) leads to a 
point in the intersection
\[
u \in \bigcap_{j=0}^{\infty} \bB_{u_j}(r_j) \ .
\]
Starting from $u$, one obtains a connecting orbit. In a Banach space setting, if 
$N=\infty$, then one can choose $r_{j+1} \leq \frac{1}{2} r_j$ for any $j$ in 
(\ref{bnt}). Choosing an arbitrary point $v_j$ in $\bB_{u_j}(r_j)$, one gets a 
Cauchy sequence $\{ v_j \}$. Thus
\[
\lim_{j\ra \infty}v_j = v \in \bigcap_{j=0}^{\infty} \bB_{u_j}(r_j) \ .
\]
Starting from $v$, one still obtains a connecting orbit.
\end{remark}

\subsection{Arnold Diffusion of DNLS ($N=3$, Non-resonant Case)}

In this subsection, we prove the existence of Arnold diffusion for a perturbed 
DNLS when $N=3$, which is a 5-dimensional system. For arbitrary $N$, one can 
find large enough annular region inside the invariant plane $\Pi$ (\ref{invp}), which 
is normally hyperbolic, for which the current proof can be easily applied. 
The point is that increasing unstable and stable dimensions does not pose 
substantial computational difficulty to establishing Arnold diffusion, while 
increasing the dimensions of tori does. We will study here the case 
that there is no resonance ($\om = 0$) inside the invariant plane $\Pi$ (\ref{invp}).
The resonant case ($\om \neq 0$) will be studied in next subsection.

\begin{figure}
\includegraphics{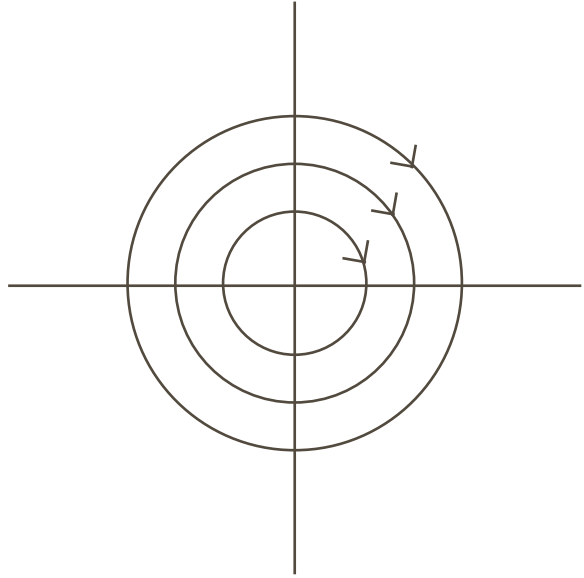}
\caption{Dynamics inside $\Pi$ (non-resonant case).}
\label{dpin}
\end{figure}

\begin{figure}
\includegraphics{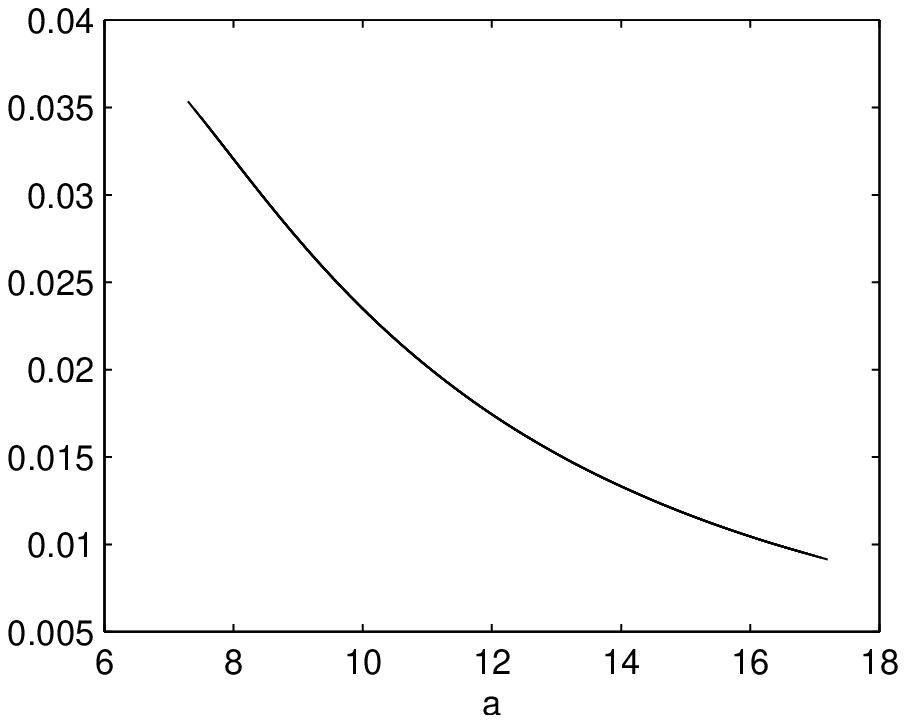}
\caption{The graph of $\sqrt{M_1^2+M_2^2}$ as a function of $a$ in the non-resonant
case $\om =0$.}
\label{nr1}
\end{figure}

\begin{figure}
\includegraphics{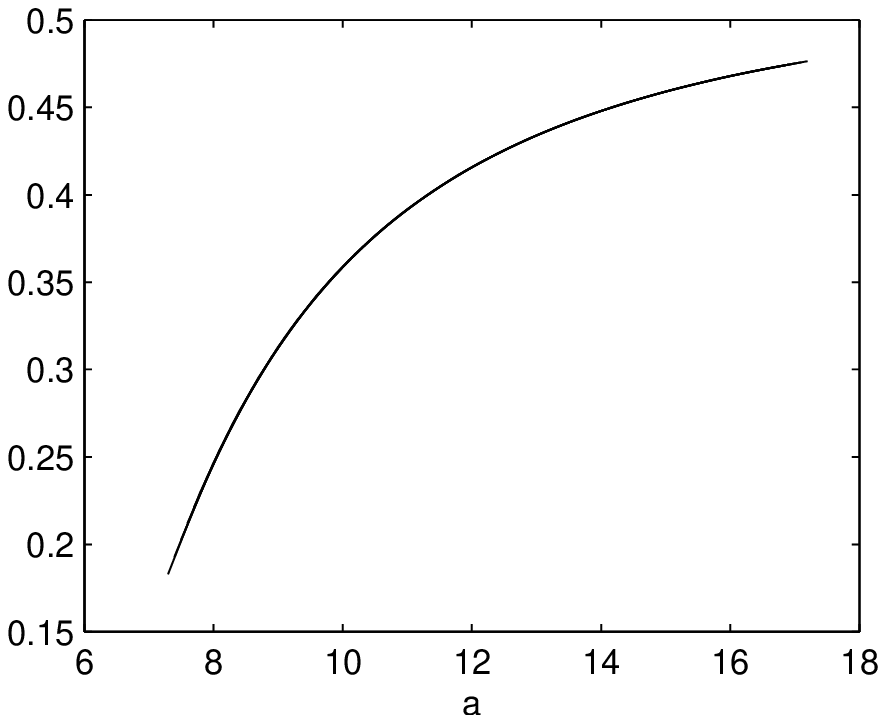}
\caption{The graph of $\sqrt{M_3^2+M_4^2}$ as a function of $a$ in the non-resonant
case $\om =0$.}
\label{nr2}
\end{figure}

\begin{figure}
\includegraphics{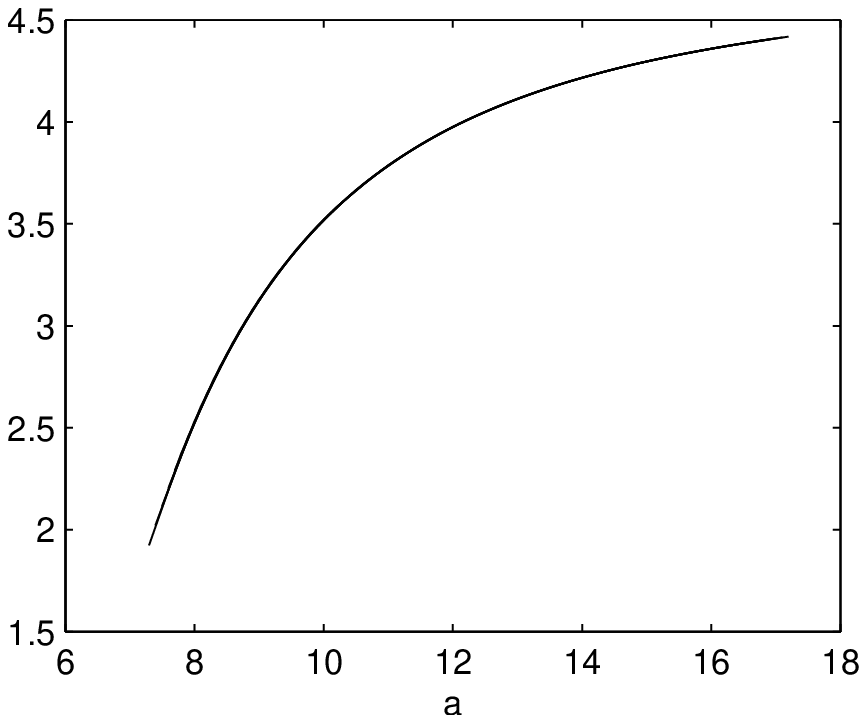}
\caption{The graph of $\sqrt{M_5^2+M_6^2}$ as a function of $a$ in the non-resonant
case $\om =0$.}
\label{nr3}
\end{figure}

Consider the following perturbation of the DNLS (\ref{DNLS})
\[
H = H_0 +\e H_1 \ , 
\]
where
\[
H_1 = \al \sin t \sum_{n=0}^{N-1} 
\left | \frac{q_n-q_{n-1}}{h}\right |^2 + \sum_{n=0}^{N-1} 
\left [ \left ( \frac{q_n-q_{n-1}}{h}\right )^2 + 
\left ( \frac{\overline{q_n}-\overline{q_{n-1}}}{h}\right )^2\right ] \ ,
\]
where $\al$ is a real parameter. Under this perturbation, dynamics inside $\Pi$ is 
unchanged. $\Pi$ consists of periodic orbits forming concentric circles (\ref{us}) 
[Figure \ref{dpin}]. 
We are interested in the following normally hyperbolic annular region inside $\Pi$
\[
\A = \left \{ \vq \in \Pi \ | \quad q_n =q, \ \forall n, \quad 3\tan{\pi \over 3}< |q| < B \right \} 
\]
where $B$ is an arbitrary large constant. Denote by $\{ \F^{u,s}(\vq): \ \vq \in \A \}$ the 
$C^1$ families of $C^2$ one dimensional unstable and stable Fenichel fibers with base points in $\A$ 
\cite{Li04} such that for any $\vq_* \in \F^u(\vq )$ or $\vq_* \in \F^s(\vq )$,
( $\vq \in \A$),
\[
\| F^t(\vq_* ) - \vq \| \leq C e^{\k t} \| \vq_* - \vq \| \ , \quad \forall t \in (-\infty , 0] \ ,
\]
or 
\[
\| F^t(\vq_* ) - \vq \| \leq C e^{-\k t} \| \vq_* - \vq \| \ , \quad \forall t \in [0, \infty ) \ ,
\]
where $F^t$ is the evolution operator of the perturbed DNLS, $\k$ and $C$ are some 
positive constants. The Fenichel fibers are $C^1$ in $\e \in [0, \e_0)$ for some $\e_0 >0$. 
It turns out that the constant of motion of DNLS (\ref{DNLS}) 
\[
I = \frac{1}{h^2} \sum_{n=0}^{N-1} \ln \rho_n 
\]
and $F_1$ [cf: (\ref{deff}) and (\ref{melv})] are 
the best choices to build the two Melnikov-Arnold intergals. Restricted to $\Pi$,
\[
I = \frac{1}{h^3} \ln \rho \ , \quad  \rho = 1+ h^2 |q|^2 \ .
\]
The level sets of $I$ lead to all the periodic orbits (1-tori) in $\Pi$. The unstable 
and stable manifolds of an 1-torus given by $I = A$ (a constant) in $\A$ are 
\[
W^{u,s}(A) = \bigcup_{\vq \in \A ,\ I(\vq ) = A} \F^{u,s}(\vq) \ ,
\]
which are three dimensional (taking into account the time dimension). 
\begin{theorem}[Arnold Diffusion]
For any $A_1$ and $A_2$ such that 
\[
\frac{1}{h^3} \ln \rho_0 < A_1 < A_2 < +\infty \ ,
\]
where $\rho_0 = 1+ h^2 \left ( 3\tan \frac{\pi}{3}\right )^2$, 
there exists a $\al_0 >0$ such that when $|\al | > \al_0$, 
$W^u(A_1)$ and $W^s(A_2)$ are connected by an orbit.
\label{nrthm}
\end{theorem}
\begin{proof}
One can check directly that for any $\vq \in \A$, $F_1 (\vq ) = -2$ and 
$\pa F_1 (\vq ) / \pa \vq = 0$ [cf: (\ref{deff}) and (\ref{melv})].
Now consider $W^s(a_1)$ and $W^u(a_2)$. Along any 
orbit $\vq^{\ s}(t)$ in $W^s(a_1)$, we have 
\begin{eqnarray*}
& & \lim_{t \ra +\infty} F_1 (\vq^{\ s} (t)) - F_1 (\vq^{\ s} (t)) = -2 - F_1 (\vq^{\ s} (t)) \\
& & = \int_{t}^{+\infty} \frac{dF_1}{dt} dt = -i\e \int_{t}^{+\infty}\{ F_1, H_1 \} dt\ ,  \\
& & \lim_{t \ra +\infty} I (\vq^{\ s} (t)) - I (\vq^{\ s} (t)) = a_1 - I (\vq^{\ s} (t)) \\
& & = \int_{t}^{+\infty} \frac{dI}{dt} dt = -i\e \int_{t}^{+\infty}\{ I, H_1 \} dt\ , 
\end{eqnarray*}
where 
\[
\{ f,g \} = \sum_{n=0}^{N-1} \rho_n \left [ \frac{\pa f}{\pa q_n} \frac{\pa g}{\pa \overline{q_n}}
- \frac{\pa f}{\pa \overline{q_n}} \frac{\pa g}{\pa q_n} \right ]
\]
is the Poisson bracket. Notice that $\{ F_1, H_0 \} = \{ I, H_0 \} = 0$ at any $\vq \in \cS$. 
Since $\pa H_1 /\pa \vq \ra 0$ exponentially as $t \ra +\infty$, the corresponding integrals converge. 
Similarly along any orbit $\vq^{\ u}(t)$ in $W^u(a_2)$, we have
\begin{eqnarray*}
& & F_1 (\vq^{\ u} (t)) - \lim_{t \ra -\infty} F_1 (\vq^{\ u} (t)) = F_1 (\vq^{\ u}(t))+2 \\
& & = \int^{t}_{-\infty} \frac{dF_1}{dt} dt = -i\e \int^{t}_{-\infty}\{ F_1, H_1 \} dt\ ,  \\
& & I (\vq^{\ u} (t)) - \lim_{t \ra -\infty} I (\vq^{\ u} (t)) = I (\vq^{\ u} (t)) - a_2 \\
& & = \int^{t}_{-\infty} \frac{dI}{dt} dt = -i\e \int^{t}_{-\infty}\{ I, H_1 \} dt\ . 
\end{eqnarray*}
Thus a neighborhood of $W^s(a_1)$ in $\cS$ can be parameterized by ($\ga , t_0 , t , F_1 , I$)
where $\ga$ is defined in (\ref{us}), $t_0$ is the initial time, and 
\[
F_1 = F_1 (\vq^{\ s} (t)) + v_1^s \ , \quad  I = I (\vq^{\ s} (t)) + v_2^s \ .
\]
When $v_1^s = v_2^s = 0$, we get $W^s(a_1)$. Thus $W^s(a_1)\cap W^u(a_2) \neq \emptyset$ 
if and only if 
\[
F_1 (\vq^{\ s} (t))=F_1 (\vq^{\ u} (t))\ , \quad I (\vq^{\ s} (t)) = I (\vq^{\ u} (t)) \ ,
\]
for some $\ga$ and $t_0$. In such a case, there is an orbit 
$\vq (t,\e ) \subset W^s(a_1)\cap W^u(a_2)$ along which 
\[
\int^{+\infty}_{-\infty}\{ F_1, H_1 \} |_ {\vq (t,\e )}dt = 0 \ , 
\quad a_1 -a_2 = -i\e \int_{-\infty}^{+\infty}\{ I, H_1 \} |_ {\vq (t,\e )}dt \ .
\]
Let $\vq (t,0 )$ be an orbit of DNLS such that $\vq (0,0 )$ and  $\vq (0,\e )$ have the 
same stable fiber base point. Then 
\[
\|  \vq (0,0 ) - \vq (0,\e ) \| \sim \O (\e ) \ .
\]
For any small $\dl >0$, there is a $T>0$ such that 
\[
\left | \int^{\pm T}_{\pm \infty}\{ F_1, H_1 \} |_ {\vq (t,\e )}dt \right | < \dl \ , 
\quad 
\left | \int^{\pm T}_{\pm \infty}\{ I, H_1 \} |_ {\vq (t,\e )}dt  \right | < \dl \ , 
\quad \forall \e \in [0,\e_0]\ ,
\]
for some $\e_0 >0$. For this $T$, when $\e$ is sufficiently small,
\[
\|  \vq (t,0 ) - \vq (t,\e ) \| \sim \O (\e ) \ , \quad \forall t \in [-T,T]\ .
\]
Thus 
\begin{eqnarray*}
& & \int^{+ T}_{-T}\{ F_1, H_1 \} |_ {\vq (t,\e )}dt =
\int^{+ T}_{-T}\{ F_1, H_1 \} |_ {\vq (t,0 )}dt + \O (\e ) \ , \\
& & \int^{+ T}_{-T}\{ I, H_1 \} |_ {\vq (t,\e )}dt =
\int^{+ T}_{-T}\{ I, H_1 \} |_ {\vq (t,0 )}dt + \O (\e ) \ .
\end{eqnarray*}
Finally we have 
\begin{eqnarray}
& & \int^{+\infty}_{-\infty}\{ F_1, H_1 \} |_ {\vq (t,\e )}dt = 
\int^{+\infty}_{-\infty}\{ F_1, H_1 \} |_ {\vq (t,0 )}dt + \O (\dl ) = 0 \ ,  \label{meq1}\\
& & a_1 -a_2 = -i\e \int_{-\infty}^{+\infty}\{ I, H_1 \} |_ {\vq (t,\e )}dt
= -i\e \int_{-\infty}^{+\infty}\{ I, H_1 \} |_ {\vq (t,0 )}dt + \O (\dl \e ) \ . \label{meq2}
\end{eqnarray}
Next we solve the above equations at the leading order in $\dl$ and $\e$. Rewrite the derivative 
given in (\ref{melv}) as follows
\[
\pa F_1 / \pa q_n = V_n e^{-i\hga}\ ,
\]
where $\hga = \ga + 2(a^2-\om^2)t_0$, $t_0 = p /\mu$, and $V_n$ represents the rest 
which does not depend on $\hga$. We also rewrite $q_c$ (\ref{us}) and $q_n$ (\ref{horbit}) 
as 
\[
q_c = \hq_c e^{i\hga}\ , \quad q_n = \hq_n e^{i\hga}\ ,
\]
where $\hq_c$ and $\hq_n$ do not depend on $\hga$. Then substitute all these into the leading 
order terms in (\ref{meq1})-(\ref{meq2}), we obtain the following equations
\begin{eqnarray}
& & \al \sqrt{M_1^2+M_2^2} \sin (t_0 +\th_1) + \sqrt{M_3^2+M_4^2} \sin (2\hga +\th_2) = 0 \ ,
\label{mele1} \\
& & a_1-a_2 = 2\e \sqrt{M_5^2+M_6^2} \sin (2\hga +\th_3) \ , \label{mele2}
\end{eqnarray}
where 
\begin{eqnarray*}
& & \cos \th_1 = \frac{M_1}{\sqrt{M_1^2+M_2^2}} \ , \quad 
\sin \th_1 = \frac{M_2}{\sqrt{M_1^2+M_2^2}} \ , \quad 
\cos \th_2 = \frac{M_3}{\sqrt{M_3^2+M_4^2}} \ ,  \\
& & \sin \th_2 = \frac{M_4}{\sqrt{M_3^2+M_4^2}} \ , \quad 
\cos \th_3 = \frac{M_5}{\sqrt{M_5^2+M_6^2}} \ , \quad 
\sin \th_3 = \frac{M_6}{\sqrt{M_5^2+M_6^2}} \ ,  \\
& & M_1 = \int_{-\infty}^{+\infty} \sum_{n=0}^{N-1} \cos \tau \ \rho_n \text{ Im } 
[V_n G_n^1 ] \ d \tau \ , \\
& & M_2 = -\int_{-\infty}^{+\infty} \sum_{n=0}^{N-1} \sin \tau \ \rho_n \text{ Im } 
[V_n G_n^1 ] \ d \tau \ , \\
& & M_3 = \int_{-\infty}^{+\infty} \sum_{n=0}^{N-1} \rho_n \text{ Re } 
[V_n G_n^2 ]\ d \tau \ , \quad 
M_4 = -\int_{-\infty}^{+\infty} \sum_{n=0}^{N-1} \rho_n \text{ Im } 
[V_n G_n^2 ] \ d \tau \ , \\
& & M_5 = -\int_{-\infty}^{+\infty} \sum_{n=0}^{N-1} \text{ Re } G_n^3 \ d \tau \ , \quad 
M_6 = \int_{-\infty}^{+\infty} \sum_{n=0}^{N-1} \text{ Im } G_n^3 \ d \tau \ , \\
& & G_n^1 = \frac{\hq_{n+1} - 2\hq_n + \hq_{n-1}}{h^2} \ , \quad 
G_n^2 = 2\frac{\overline{\hq_{n+1}} - 2\overline{\hq_n} + \overline{\hq_{n-1}}}{h^2} \ , \\
& & G_n^3 = -2\frac{\overline{\hq_n}(\overline{\hq_{n+1}} - 2\overline{\hq_n} + 
\overline{\hq_{n-1}})}{h^2} \ ,
\end{eqnarray*}
and $\tau = t +p/\mu$. Equations (\ref{mele1})-(\ref{mele2}) are easily solvable as long as
neither $\sqrt{M_1^2+M_2^2}$ nor $\sqrt{M_5^2+M_6^2}$ vanishes. In Figures \ref{nr1}-\ref{nr3}, we 
plot the graphs of them as functions of $a$. We solve equation (\ref{mele2}) for $\hga$, 
then solve equation (\ref{mele1}) for $t_0$. Thus when $|\al |$ is large enough, we have 
solutions. It is also clear from equations (\ref{mele1})-(\ref{mele2}) that $W^s(a_1)$ and 
$W^u(a_2)$ intersect transversally. Then we can choose a sequence
\[
A_1=a_1 < a_2 < \cdots < a_N=A_2 \ ,
\]
such that $W^s(a_j)$ and $W^u(a_{j+1})$ ($ 1\leq j \leq N-1$) intersect transversally.
The period of the 1-tori (\ref{us}) is $\pi / a^2$. Thus we can always choose the $a_j$'s 
such that the frequencies $\frac{1}{2a^2}$ of the corresponding 1-tori are irrational. 
Therefore we obtain a transition chain. Apply Lemma \ref{al} to the period-$2\pi$ map 
of the DNLS, we obtain the claim of the theorem.
\end{proof}
\begin{remark}
The constant of motion $I$ is equivalent to $F_2$ for $z_2^c =1$ (\ref{deff}).
The continuum limit of $H_1$ has the form 
\[
H_1 = \al \sin t \int_0^1 |q_x|^2 dx + \int_0^1 (q_x^2 + \overline{q_x}^{\ 2})dx \ ,
\]
which is suitable for the NLS setting. One can regularize the perturbation by replacing the 
partial derivative $\pa_x$ in $H_1$ by a Fourier multiplier $\hat{\pa}_x$, e.g. a 
Galerkin truncation. One can use the constant of motion 
\[
I = \int_0^1 |q|^2 dx 
\]
to build the second Melnikov-Arnold integral. 
\end{remark}


\begin{figure}
\includegraphics{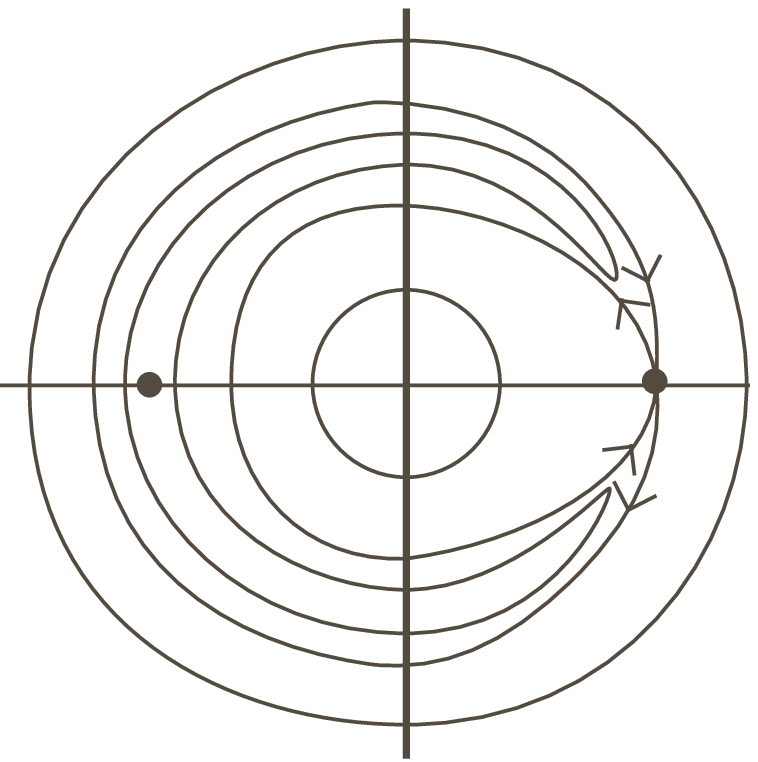}
\caption{Dynamics inside $\Pi$ (resonant case).}
\label{dpir}
\end{figure}

\begin{figure}
\includegraphics{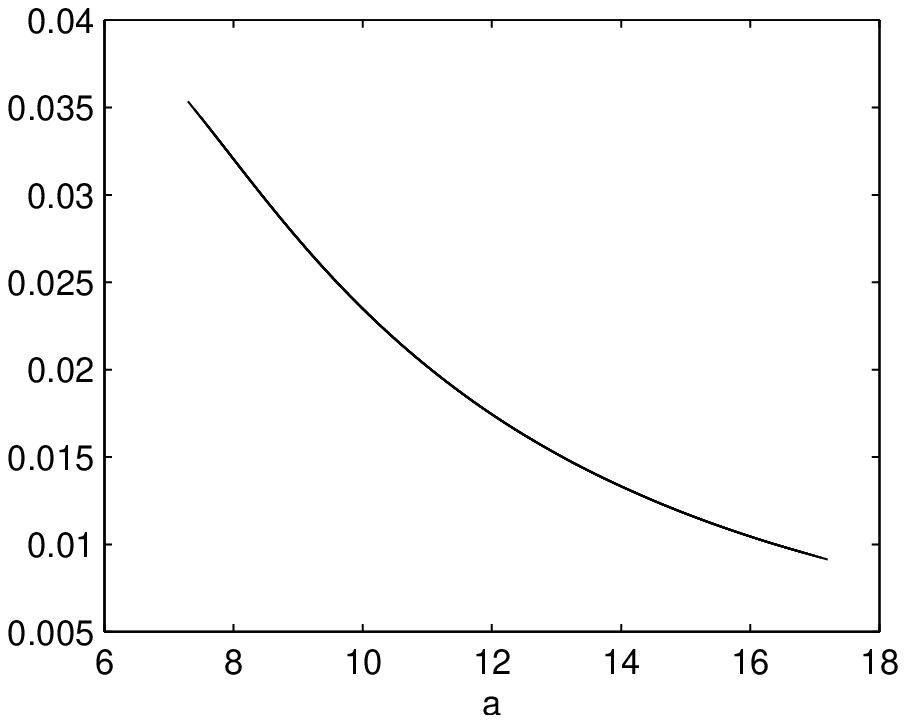}
\caption{The graph of $\sqrt{M_1^2+M_2^2}$ as a function of $a$ in the resonant
case $\om =10$.}
\label{r1}
\end{figure}
\begin{figure}
\includegraphics{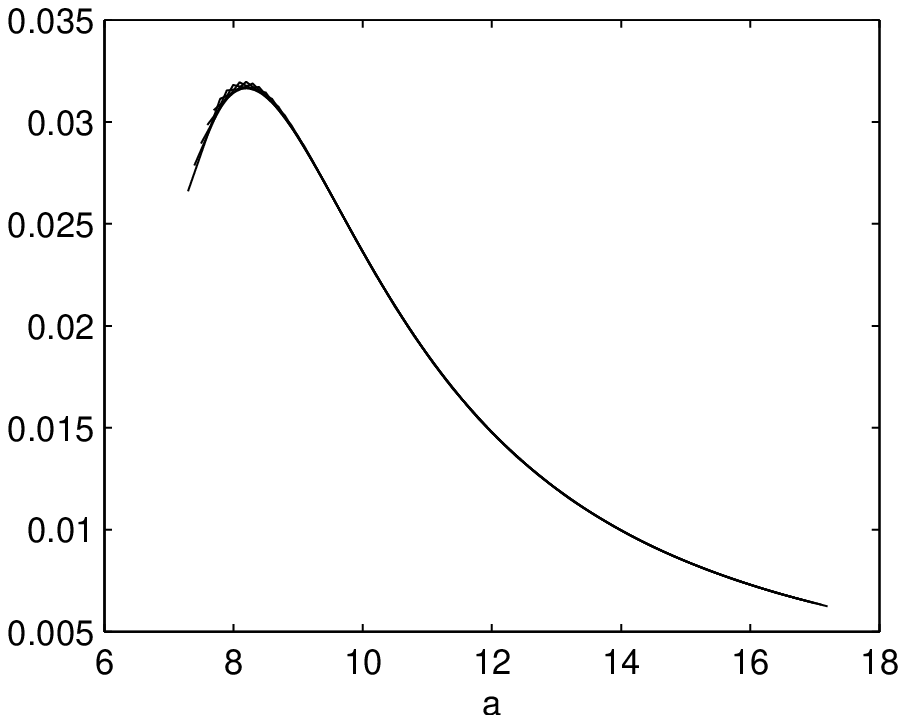}
\caption{The graph of $\sqrt{M_3^2+M_4^2}$ as a function of $a$ in the resonant
case $\om =10$.}
\label{r2}
\end{figure}
\begin{figure}
\includegraphics{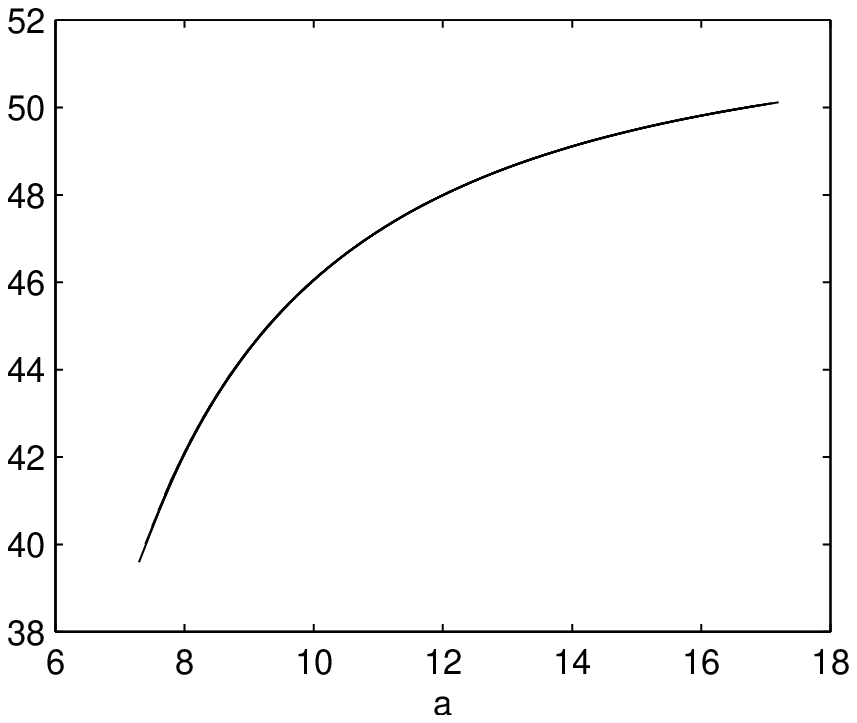}
\caption{The graph of $\sqrt{M_5^2+M_6^2}$ as a function of $a$ in the resonant
case $\om =10$.}
\label{r3}
\end{figure}

\subsection{Arnold Diffusion of DNLS ($N=3$, Resonant Case)}

In this subsection, we prove the existence of Arnold diffusion for a perturbed 
DNLS when $N=3$, which is a 5-dimensional system. We will study here the case 
that there is a resonance ($\om > 3\tan{\pi \over 3}$) inside the invariant 
plane $\Pi$ (\ref{invp}).

Consider the following perturbation of the DNLS (\ref{DNLS})
\[
H = H_0 +\e (H_1 +H_2) \ , 
\]
where
\[
H_1 = \al \sum_{n=0}^{N-1}(q_n +\overline{q_n}) \ , \quad 
H_2 = \sin t \sum_{n=0}^{N-1}
\left | \frac{q_n-q_{n-1}}{h}\right |^2 \ ,
\]
where $\al$ is a real parameter. Under this perturbation, dynamics inside $\Pi$ is 
changed. Due to the resonance $a =\om$ in (\ref{us}), some tori do not persist into 
KAM tori. A secondary separatrix is generated. Inside this separatrix are the 
secondary tori [Figure \ref{dpir}]. As can be seen below, resonance does not add 
difficulty to the Arnold diffusion problem. Instead of $I$ in last subsection, we use 
\[
\hH = H = H_0 +\e H_1 \ ,
\]
to build one of the two Melnikov-Arnold integrals. 
Restricted to $\Pi$, the level sets of $\hH$ produces Figure \ref{dpir}. The unstable 
and stable manifolds of an 1-torus given by $\hH = A$ (a constant) in $\A$ are 
\[
W^{u,s}(A) = \bigcup_{\vq \in \A ,\ \hH(\vq ) = A} \F^{u,s}(\vq) \ .
\]
Let 
\[
A_* = \frac{1}{h^3}\left [ 2\left ( 3\tan \frac{\pi}{3}\right )^2
-\frac{2}{h^2}(1+\om^2h^2) \ln \rho_0 \right ]\ ,
\]
where $\rho_0 = 1+ h^2 \left ( 3\tan \frac{\pi}{3}\right )^2$.
\begin{theorem}[Arnold Diffusion]
For any $A_1$ and $A_2$ such that 
\[
A_* < A_1 < A_2 < +\infty \ ,
\]
there exists a $\al_0 >0$ such that when $|\al | > \al_0$, 
$W^u(A_1)$ and $W^s(A_2)$ are connected by an orbit.
\end{theorem}
\begin{proof}
Again one can check directly that for any $\vq \in \A$, $F_1 (\vq ) = -2$ and 
$\pa F_1 (\vq ) / \pa \vq = 0$ [cf: (\ref{deff}) and (\ref{melv})].
Similar to the proof of Theorem \ref{nrthm}, consider $W^s(a_1)$ and $W^u(a_2)$. 
Along any orbit $\vq^{\ s}(t)$ in $W^s(a_1)$, we have 
\begin{eqnarray*}
& & \lim_{t \ra +\infty} F_1 (\vq^{\ s} (t)) - F_1 (\vq^{\ s} (t)) = -2 - F_1 (\vq^{\ s} (t)) \\
& & = \int_{t}^{+\infty} \frac{dF_1}{dt} dt = -i\e \int_{t}^{+\infty}\{ F_1, H_1+H_2 \} dt\ ,  \\
& & \lim_{t \ra +\infty} \hH (\vq^{\ s} (t)) - \hH (\vq^{\ s} (t)) = a_1 - \hH (\vq^{\ s} (t)) \\
& & = \int_{t}^{+\infty} \frac{d\hH}{dt} dt = -i\e \int_{t}^{+\infty}\{ \hH, H_2 \} dt\ , 
\end{eqnarray*}
where 
\[
\{ f,g \} = \sum_{n=0}^{N-1} \rho_n \left [ \frac{\pa f}{\pa q_n} \frac{\pa g}{\pa \overline{q_n}}
- \frac{\pa f}{\pa \overline{q_n}} \frac{\pa g}{\pa q_n} \right ]
\]
is the Poisson bracket. Notice that $\{ F_1, H_0 \} = \{ \hH , \hH\} = 0$ at any $\vq \in \cS$. 
Since $\frac{\pa F_1}{\pa \vq}, \frac{\pa H_2}{\pa \vq} \ra 0$ exponentially as $t \ra +\infty$, 
the corresponding integrals converge. 
Similarly along any orbit $\vq^{\ u}(t)$ in $W^u(a_2)$, we have
\begin{eqnarray*}
& & F_1 (\vq^{\ u} (t)) - \lim_{t \ra -\infty} F_1 (\vq^{\ u} (t)) = F_1 (\vq^{\ u}(t))+2 \\
& & = \int^{t}_{-\infty} \frac{dF_1}{dt} dt = -i\e \int^{t}_{-\infty}\{ F_1, H_1+H_2 \} dt\ ,  \\
& & \hH (\vq^{\ u} (t)) - \lim_{t \ra -\infty} \hH (\vq^{\ u} (t)) = \hH (\vq^{\ u} (t)) - a_2 \\
& & = \int^{t}_{-\infty} \frac{d\hH}{dt} dt = -i\e \int^{t}_{-\infty}\{ \hH, H_2 \} dt\ . 
\end{eqnarray*}
Thus a neighborhood of $W^s(a_1)$ in $\cS$ can be parameterized by ($\vth , t_0 , t , F_1 , \hH$)
where $\vth$ is the angle of the 1-torus $\hH = a_1$ in $\Pi$, $t_0$ is the initial time, and 
\[
F_1 = F_1 (\vq^{\ s} (t)) + v_1^s \ , \quad  \hH = \hH (\vq^{\ s} (t)) + v_2^s \ .
\]
When $v_1^s = v_2^s = 0$, we get $W^s(a_1)$. Thus $W^s(a_1)\cap W^u(a_2) \neq \emptyset$ 
if and only if 
\[
F_1 (\vq^{\ s} (t))=F_1 (\vq^{\ u} (t))\ , \quad  \hH (\vq^{\ s} (t)) = \hH (\vq^{\ u} (t)) \ ,
\]
for some $\vth$ and $t_0$. In such a case, there is an orbit 
$\vq (t,\e ) \subset W^s(a_1)\cap W^u(a_2)$ along which 
\[
\int^{+\infty}_{-\infty}\{ F_1, H_1+H_2 \} |_ {\vq (t,\e )}dt = 0 \ , 
\quad a_1 -a_2 = -i\e \int_{-\infty}^{+\infty}\{ \hH , H_2 \} |_ {\vq (t,\e )}dt \ .
\]
Let $\vq (t,0 )$ be an orbit of DNLS such that $\vq (0,0 )$ and  $\vq (0,\e )$ have the 
same stable fiber base point. Then 
\[
\|  \vq (0,0 ) - \vq (0,\e ) \| \sim \O (\e ) \ .
\]
For any small $\dl >0$, there is a $T>0$ such that 
\[
\left | \int^{\pm T}_{\pm \infty}\{ F_1, H_1+H_2 \} |_ {\vq (t,\e )}dt \right | < \dl \ , 
\quad 
\left | \int^{\pm T}_{\pm \infty}\{ \hH , H_2 \} |_ {\vq (t,\e )}dt  \right | < \dl \ , 
\quad \forall \e \in [0,\e_0]\ ,
\]
for some $\e_0 >0$. For this $T$, when $\e$ is sufficiently small,
\[
\|  \vq (t,0 ) - \vq (t,\e ) \| \sim \O (\e ) \ , \quad \forall t \in [-T,T]\ .
\]
Thus 
\begin{eqnarray*}
& & \int^{+ T}_{-T}\{ F_1, H_1+H_2 \} |_ {\vq (t,\e )}dt =
\int^{+ T}_{-T}\{ F_1, H_1+H_2 \} |_ {\vq (t,0 )}dt + \O (\e ) \ , \\
& & \int^{+ T}_{-T}\{ \hH , H_2 \} |_ {\vq (t,\e )}dt =
\int^{+ T}_{-T}\{ \hH , H_2 \} |_ {\vq (t,0 )}dt + \O (\e ) \ .
\end{eqnarray*}
Finally we have 
\begin{eqnarray}
& & \int^{+\infty}_{-\infty}\{ F_1, H_1+H_2 \} |_ {\vq (t,\e )}dt = 
\int^{+\infty}_{-\infty}\{ F_1, H_1+H_2 \} |_ {\vq (t,0 )}dt + \O (\dl ) = 0 \ ,  \label{mer1}\\
& & a_1 -a_2 = -i\e \int_{-\infty}^{+\infty}\{ \hH , H_2 \} |_ {\vq (t,\e )}dt
= -i\e \int_{-\infty}^{+\infty}\{ H_0, H_2 \} |_ {\vq (t,0 )}dt + \O (\dl \e ) \ . \label{mer2}
\end{eqnarray}
To the leading order terms in (\ref{mer1})-(\ref{mer2}), we obtain the following equations
\begin{eqnarray}
& & \sqrt{M_1^2+M_2^2} \sin (t_0 +\th_1) + \al \sqrt{M_3^2+M_4^2} \sin (\hga +\th_2) = 0 \ ,
\label{melr1} \\
& & a_1-a_2 = 2\e \sqrt{M_5^2+M_6^2} \sin (t_0 +\th_3) \ , \label{melr2}
\end{eqnarray}
where 
\begin{eqnarray*}
& & \cos \th_1 = \frac{M_1}{\sqrt{M_1^2+M_2^2}} \ , \quad 
\sin \th_1 = \frac{M_2}{\sqrt{M_1^2+M_2^2}} \ , \quad 
\cos \th_2 = \frac{M_3}{\sqrt{M_3^2+M_4^2}} \ ,  \\
& & \sin \th_2 = \frac{M_4}{\sqrt{M_3^2+M_4^2}} \ , \quad 
\cos \th_3 = \frac{M_5}{\sqrt{M_5^2+M_6^2}} \ , \quad 
\sin \th_3 = \frac{M_6}{\sqrt{M_5^2+M_6^2}} \ ,  \\
& & M_1 = \int_{-\infty}^{+\infty} \sum_{n=0}^{N-1} \cos \tau \ \rho_n \text{ Im } 
[V_n G_n^1 ] \ d \tau \ , \\
& & M_2 = -\int_{-\infty}^{+\infty} \sum_{n=0}^{N-1} \sin \tau \ \rho_n \text{ Im } 
[V_n G_n^1 ] \ d \tau \ , \\
& & M_3 = -\int_{-\infty}^{+\infty} \sum_{n=0}^{N-1} \rho_n \text{ Re } 
[V_n]\ d \tau \ , \quad 
M_4 = \int_{-\infty}^{+\infty} \sum_{n=0}^{N-1} \rho_n \text{ Im } 
[V_n] \ d \tau \ , \\
& & M_5 = \int_{-\infty}^{+\infty} \sum_{n=0}^{N-1} \cos \tau \text{ Im } [G_n^1G_n^2] \ d \tau \ , \\
& & M_6 = -\int_{-\infty}^{+\infty} \sum_{n=0}^{N-1} \sin \tau \text{ Im } [G_n^1G_n^2] \ d \tau \ , \\
& & G_n^1 = \frac{\hq_{n+1} - 2\hq_n + \hq_{n-1}}{h^2} \ , \\
& & G_n^2 = \frac{1}{h^2}\left [\overline{\hq_{n+1}} - 2\overline{\hq_n} + \overline{\hq_{n-1}}
\right ] 
+ |\hq_n|^2 \left [\overline{\hq_{n+1}} + \overline{\hq_{n-1}}\right ] -2\om^2 \overline{\hq_n}\ ,
\end{eqnarray*}
and $\tau = t +p/\mu$. Equations (\ref{melr1})-(\ref{melr2}) are easily solvable as long as
neither $\sqrt{M_3^2+M_4^2}$ nor $\sqrt{M_5^2+M_6^2}$ vanishes. In Figures \ref{r1}-\ref{r3}, we 
plot the graphs of them as functions of $a$. We solve equation (\ref{melr2}) for $t_0$, 
then solve equation (\ref{melr1}) for $\hga$. Thus when $|\al |$ is large enough, we have 
solutions. It is also clear from equations (\ref{melr1})-(\ref{melr2}) that $W^s(a_1)$ and 
$W^u(a_2)$ intersect transversally. Then we can choose a sequence
\[
A_1=a_1 < a_2 < \cdots < a_N=A_2 \ ,
\]
such that $W^s(a_j)$ and $W^u(a_{j+1})$ ($ 1\leq j \leq N-1$) intersect transversally.
The period of the 1-tori ($\hH = a_j$) depends on $a_j$ no matter they are KAM tori or 
secondary tori. We can always choose the $a_j$'s 
such that the frequencies of the corresponding 1-tori are irrational. 
We can use one secondary torus inside and close enough to the separatrix (Figure \ref{dpir})
to bridge across the resonant region of width $\O(\sqrt{\e})$.
Therefore we obtain a transition chain. Apply Lemma \ref{al} to the period-$2\pi$ map 
of the DNLS, we obtain the claim of the theorem.
\end{proof}
\begin{remark}
The continuum limit of $H_1$ and $H_2$ have the form 
\[
H_1 = \al \int_0^1 (q +\bq ) dx \ , \quad 
H_2 = \sin t \int_0^1 |q_x|^2 dx \ ,
\]
which is suitable for the NLS setting. One can regularize the perturbation by replacing the 
partial derivative $\pa_x$ in $H_2$ by a Fourier multiplier $\hat{\pa}_x$, e.g. a 
Galerkin truncation. 
\end{remark}

\end{document}